\documentclass[10pt]{article}

\usepackage{amsmath}
\usepackage{amssymb,float}

\usepackage{graphicx}

\usepackage{cite}

\usepackage{color}

\usepackage[labelfont=bf,labelsep=period,justification=raggedright]{caption}

\makeatletter
\renewcommand{\@biblabel}[1]{\quad#1.}
\makeatother

\date{}

\newcommand{\PP}{P^{-1}}
\newcommand{\bb}{{\frac{\beta^2}{2}}}

\begin{document}

% Title must be 150 characters or less
\begin{flushleft}
{\Large \textbf{%General One-Dimensional Fokker-Planck
Complexity Reduction of
Rate-Equations Models for Two-Choice Decision-Making} }
% Insert Author names, affiliations and corresponding author email.
\\
Jos\'e Antonio Carrillo$^{1}$, St\'ephane Cordier$^{2}$, Gustavo
Deco$^{3,4}$, Simona Mancini$^{2,\ast}$
\\
\bf{1} Department of Mathematics, Imperial College London, London
SW7 2AZ, UK
\\
\bf{2} Univ. Orl\'eans, CNRS, MAPMO, UMR 7349  $\&$ FR 2964, F-45067 Orl\'eans, France
\\
\bf{3} ICREA (Instituci\'o Catalana de Recerca i Estudis Avan\c
cats), Barcelona, Spain
\\
\bf{4} Computational Neuroscience Unit, Universitat Pompeu Fabra,
Barcelona, Spain
\\
$\ast$ E-mail: simona.mancini@univ-orleans.fr
\end{flushleft}

% Please keep the abstract between 250 and 300 words
\section*{Abstract}
We are concerned with the complexity reduction of a stochastic
system of differential equations governing the dynamics of a
neuronal circuit describing a decision-making task. This reduction
is based on the slow-fast behavior of the problem and holds on the
whole phase space and not only locally around the spontaneous
state. Macroscopic quantities, such as performance and reaction
times, computed applying this reduction are in agreement with
previous works in which the complexity reduction is locally
performed at the spontaneous point by means of a Taylor expansion.

% Please keep the Author Summary between 150 and 200 words
% Use first person. PLoS ONE authors please skip this step.
% Author Summary not valid for PLoS ONE submissions.
%\section*{Author Summary}
\section*{Author Summary}
In recent years, the neuronal mechanisms governing binary
decision-making have largely been studied. In a wide range of
experiments, the behavioral response can be described by means of
stochastic systems. In particular,  we are interested in the
qualitative analysis of the performance and reaction times
corresponding to a decision-making problem. In order to perform
the link between neuronal circuits and stochastic systems, we
consider here a mean-field reduction consisting of two interacting
rate models. This system of equations can be associated with a
partial differential equation governing the probability
distribution of the activities of the different neuronal
populations. However, the nonlinear nature of the original
equations hinders analytical progress in the partial differential
framework and reduces its study to numerical investigations having
very large computational times. We propose here a complexity
reduction of the stochastic system, based on its slow-fast
characteristic, allowing us to define a solution on the whole
biological domain, rather than just on a neighborhood of the
spontaneous state, as it is usually done. Moreover, close to the
bifurcation point, the qualitative behavior of the  performance
and reaction times is the same as  the one  observed in previous
works.

\section*{Introduction}

%%%%%%%%%%%%%%%%%%%%%%%%%

Goal-oriented behavior involves the constant making of decisions
between alternative choices. Hence, a sorrow understanding of the
mechanisms underlying decision-making is fundamental in behavioral
neurosciences. During the last years tremendous advances in
neurophysiological and theoretical neurosciences have started to
reveal the neural mechanisms underlying decision-making. Most of
these studies are focused on binary decision-making. In binary
decision tasks, the subject is asked to make a choice between two
alternatives according to an experimentally defined criterion
based on sensory input information. The difficulty of the making
of the decision can be manipulated by the level of uncertainty in
the discrimination of the sensory input information which can be
influenced by regulating the signal to noise ratio. For the same
pattern of stimulation the decision of the subject from trial-to-trial
is stochastic but such that its trial-to-trial average is
determined by the input.

In a wide range of such decision-making experiments, the
behavioral response, i.e. performance and reaction times, can be
properly described by a simple stochastic diffusion model
\cite{RS04}. In the diffusion model evidences about one or the
other choice is accumulated continuously over time until a
decision boundary is reached. Hence, it is plausible to think that
the underlying neuronal system perform decision-making by
accumulation of evidences. Indeed, recent electrophysiological
recordings from awake behaving monkeys performing decision-making
have shown that trial-averaged spiking activity of some neurons
shows a gradual, nearly linear increase shortly after stimulus
onset until just before the response is made
\cite{RS01,RS03,SN01,RS02,HS05}. Let us consider for example the
so-called random-dot task. The random-dot task is a two-choice
visual decision-making task. In this task, the subject is
confronted with a visual stimulus consisting of a field of
randomly moving dots. On any one trial, a fixed fraction of the
dots, determined by the experimentalist, moves coherently in one
of two directions. The subject must discriminate in which
direction the majority of dots are moving. Electrophysiological
recordings from awake behaving monkeys performing this task have
shown that trial-averaged spiking activity of neurons in the
lateral intraparietal cortex (LIP) reflects the accumulation of
information mentioned above \cite{SN01,RS02,HS05}. In fact, LIP
cells receive excitatory drive from direction-selective cells of
extrastriate visual area MT. The difference in the firing rates of
MT cells with opposing preferred directions is a likely measure of
the available visual evidence for the direction of coherent motion
of the moving dot stimulus, further indicating that the LIP area
is involved in this decision-making task.

A computational model of the cortical circuit putatively
underlying LIP cell activity has been proposed and studied
numerically \cite{W02,WW06,LW06}. This model was able to explain
qualitatively the experimentally observed trial-averaged spiking
activity. The model consists of two populations of spiking
neurons, within which interactions are mediated by excitatory
synapses and between which interactions occur principally through
an intermediary population of inhibitory interneurons. Sensory
input may bias the competition in favor of one of the populations,
potentially resulting in a gradually developing decision in which
neurons in the chosen population exhibit increased activity while
activity in the other population is inhibited. When the activity
of one of the two populations exceeds a pre-defined threshold then
a behavioral response is generated. There are two possible
mechanisms by which the decision can be made which correspond to
two different dynamical working points. In one mechanism, the
appearance of the stimulus destabilizes the so-called spontaneous
state, in which all neurons show low activity. In the other
mechanism, even after stimulus presentation, the network can still
sustain low activity, but random fluctuations eventually drive the
collective activity of the network to an activated state, in which
a fraction of the neurons fires at a high rate. Let us remark that
the behavioral performance of the model (i.e. the fraction of
trials in which the 'correct' direction is chosen) and
corresponding reaction time to make a decision match well the
behavioral data from both monkeys and humans performing the
moving-dot task \cite{RL}. Therefore, we conjecture, that it may
be possible to establish a direct link between the neuronal and
behavioral correlates of decision-making. Here, we aim to
establish this link between the underlying physiology and the
observed behavioral response in decision-making tasks by
performing a 1D reduction of the dynamics of the neuronal circuit
in such a way that we derive from the underlying detailed neuronal
dynamics a valid "nonlinear" diffusion model.

In order to perform the explicit link between a neuronal circuit
processing decision-making and a diffusion model able to describe
the behavioral level, we consider a reduced model of the
competitive cortical circuit solving the decision-making problem.
In particular, we consider a mean-field reduction consisting of
two competing rate models. This system of Langevin equations can
be associated with a Fokker-Planck equation for the probability
distribution of the activities of the different neuronal
populations. The nonlinear nature of the original equations,
however, hinders analytical progress in the Fokker-Planck
framework. For this reason, the main analysis of such noise driven
probabilistic decision-making systems remains based on numerical
investigations, which are time consuming because of the need for
sufficiently many trials to generate statistically meaningful
data.

One dimensional Fokker-Planck effective reductions near
bifurcation points were obtained by Roxin and Ledberg \cite{RL}.
They propose to locally approximate the evolution on the slow
manifold at the bifurcation point by Taylor expanding the
nonlinearities with respect to the slow variable based on the
different time scales. In this way, they obtained an effective
potential of degree 4 for the one dimensional Fokker-Planck
dynamics on the slow manifold. We propose to build on this
strategy by not performing any Taylor expansion and rather
approximate the slow manifold in terms of the nonlinearities
defining the evolution using the slow-fast character of the
dynamical system, see \cite{BG}. Once, the slow manifold is
approximated then we project the Fokker-Planck dynamics on it to
obtain a one dimensional Fokker-Planck reduction possibly valid
beyond the local character of the expansion in \cite{RL}. More
importantly, we recover an approximated full reduced potential and
an explicit formula for the stationary state distribution on the
slow manifold. This effective reduced potential covers the other
two stable equilibrium points, and not only the central
equilibrium point changing its stability character as in
\cite{RL}. For instance, in the case of the subcritical Hopf
bifurcation in Figure \ref{bifurcation}, we get an approximated
potential as in Figure \ref{pot1D} while the approach in \cite{RL}
leads to a potential of degree 4 approximating an interval around
the central equilibrium point, see \cite[Figure 2]{RL}. Our
strategy is valid as long as the approximated slow manifold lies
fully in the relevant biological range of positive rates, this is
the reason we restrict the range of the bifurcation parameter
$w_+$ as explained in Figure \ref{manifold}.

%%%%%%%%%%%%%%% Other References: \cite{Deco,RBW,WC}

% Results and Discussion can be combined.
\section*{Results}

We showcase our strategy in a particular example of a neuronal
population model with two pools with self-excitation and
cross-inhibition as in \cite{Deco}. The firing rates $\nu_1$ and
$\nu_2$ of the neuronal networks are determined by the stochastic
dynamical system:
\begin{equation}\label{ODE}
\begin{cases}
{d\nu_1} = \left[-\nu_1 +\phi(\lambda_1+ w_+ \nu_1 + w_I
\nu_2)\right]{dt} +
\beta \,dW_t^1\\[3mm]
{d\nu_2} = \left[- \nu_2 + \phi(\lambda_2+ w_I \nu_1 + w_+
\nu_2)\right]{dt} + \beta \,dW_t^2
\end{cases},
\end{equation}
with $\nu_1, \nu_2\geq 0$. Here the applied stimuli are,
$\lambda_1 =33$, $\lambda_2=\lambda_1-\Delta \lambda$, with the
bias $\Delta \lambda\in [0, 10^{-3}]$,  the inhibitory
connectivity coefficient is $w_I=1.9$, the standard deviation is
$\beta=3\cdot 10^{-3}$, and the excitatory connectivity
coefficient $w_+$ is the bifurcation parameter, its range of
values will be discussed later on. Moreover, the sigmoid
(response) function $\phi(z)$ is given by:
$$
\phi(z)=\frac{\nu_c}{1+\exp(-b z+\alpha)}\,,
$$
with $\nu_c= 15$, $b= 0.25$, and $\alpha=11.1$.

\begin{figure}[H]
\centering{
\includegraphics[width=4in]{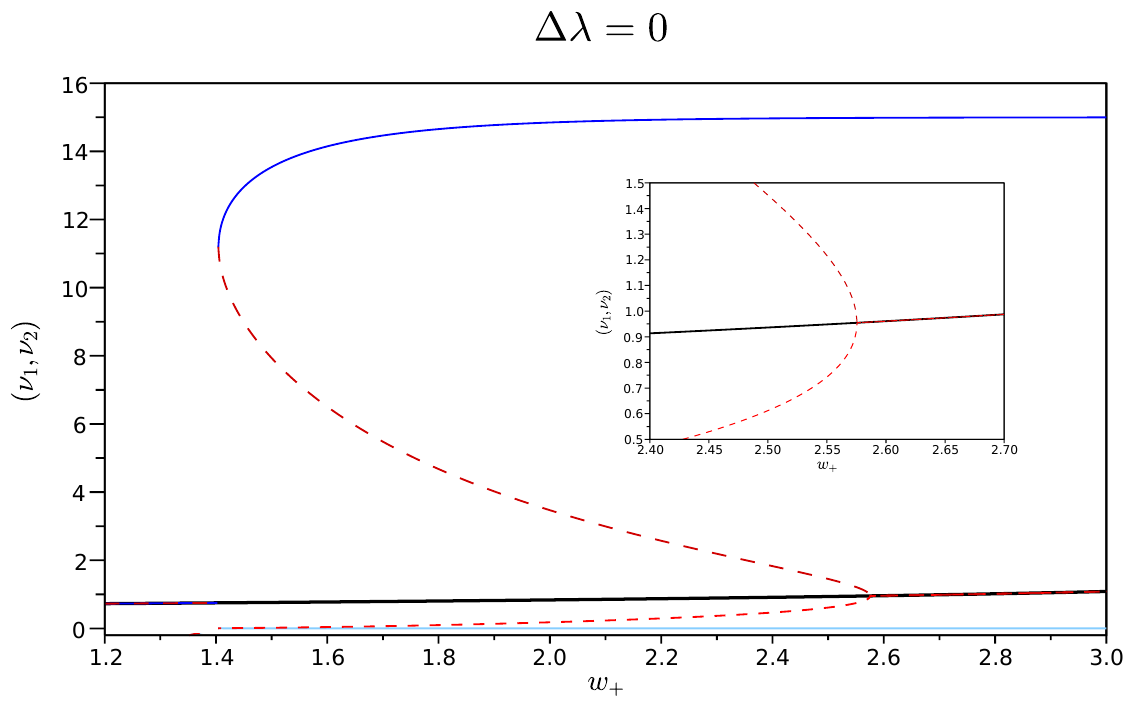}\\
\includegraphics[width=4in]{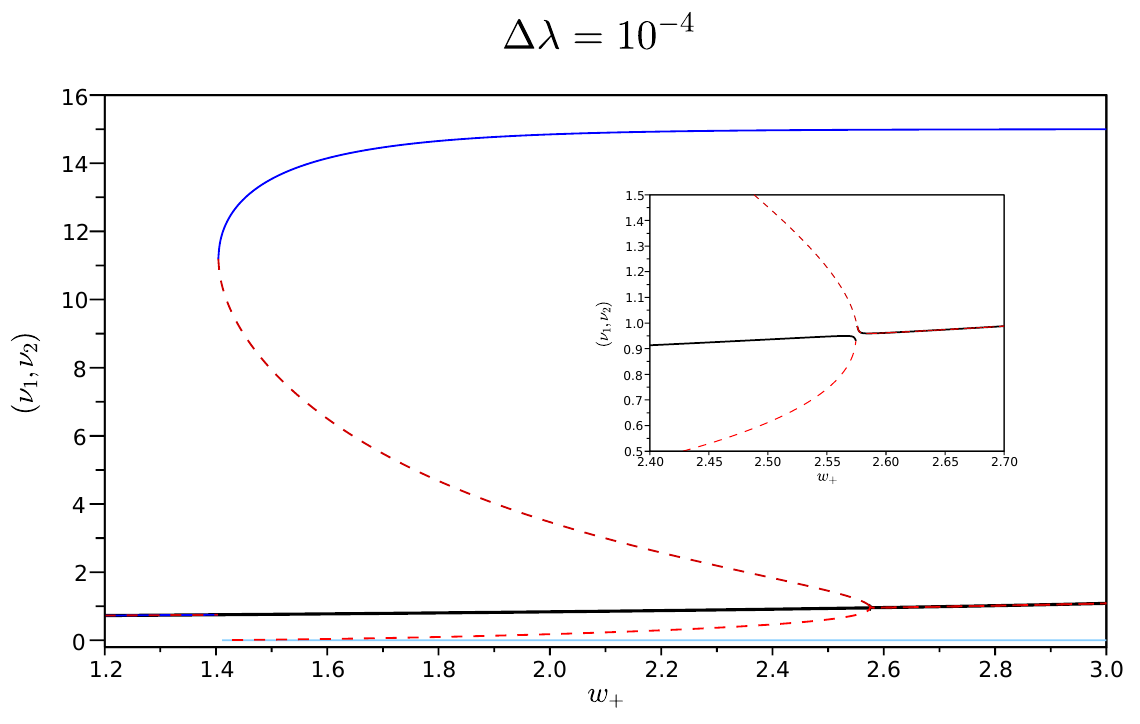}\\
\includegraphics[width=4in]{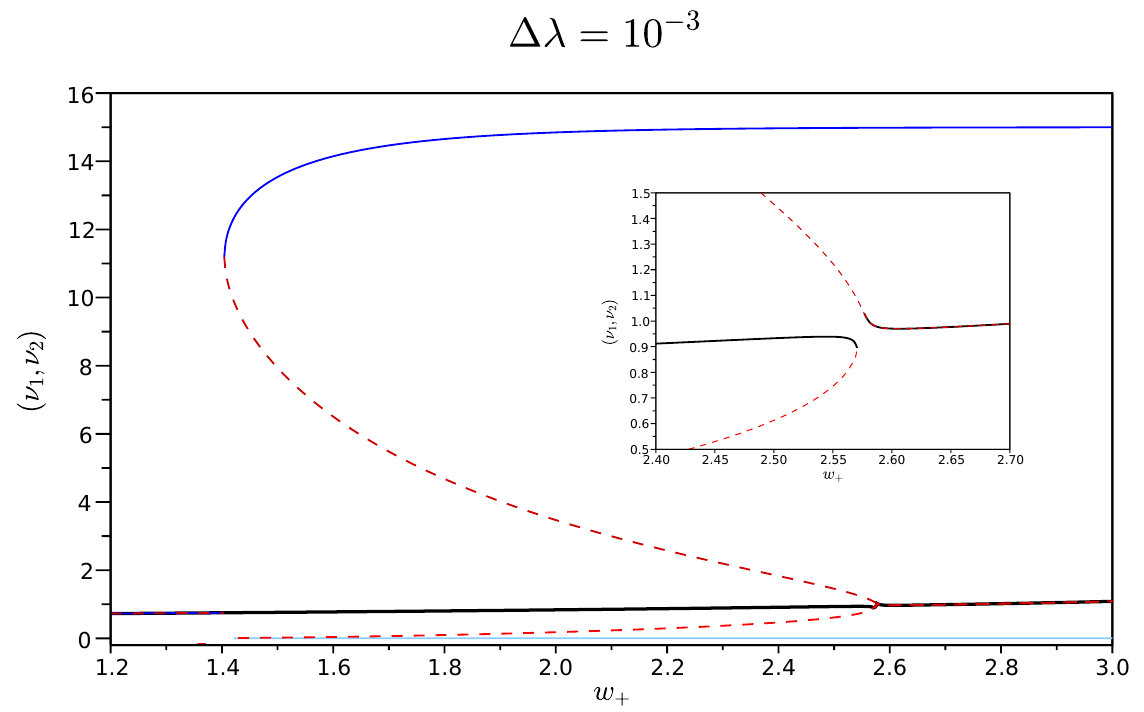}}
\caption{{\bf Bifurcation diagrams with respect to $w_+$.} From top to bottom the bias $\Delta \lambda$ takes the following values $0, 10^{-4}$ and $10^{-3}$. The components of the spontaneous state are traced by the black line. There are other two stable equilibria points which coordinates are represented by the blue and light-blue lines. The red continuous and dashed lines correspond to the unstable equilibrium points. The small picture is a zoom on the bifurcation point.} \label{bifurcation}
\end{figure}

In terms of $w_+$ the underlying dynamical system without noise
presents a subcritical Hopf bifurcation whose bifurcation diagram
is shown in Figure \ref{bifurcation}. For values of the excitatory
coefficient around $w_+=1.4$ the system pass from a single
asymptotically stable equilibria to a situation in which there are
three stable (continuous lines) and two unstable equilibria
(hashed lines). The second local bifurcation, where the central
asymptotic equilibria disappears, happens at around $w_+=2.5695$
for $\Delta \lambda=10^{-3}$.

\begin{figure}[H]
\centering{\includegraphics[width=3in]{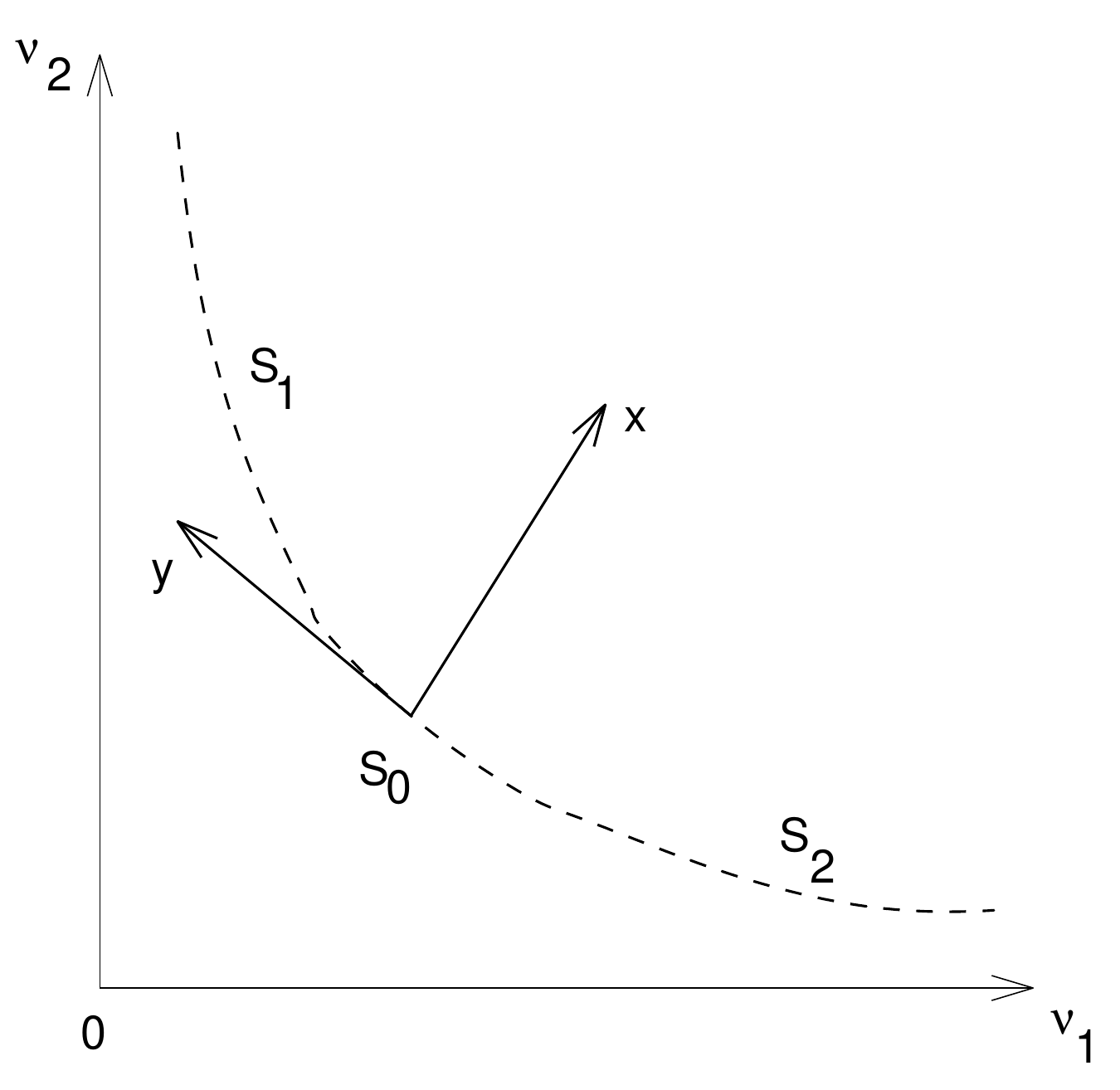}}
\caption{{\bf Change of variables.} Position of the new coordinate system $(x,y)$ with respect to the initial one $(\nu_1, \nu_2)$. The origin is taken at the central equilibrium point $S_0$ which may stable or unstable accordingly to the value of $w_+$. $S_1$ and $S_2$ are two others equilibrium points. The dashed line represents the slow manifold on which the reduced dynamics takes place, and it contains all the equilibrium points of the system. } \label{chang-variable}
\end{figure}

We next need to find an approximation of the slow manifold joining
the spontaneous $S_0$ and the decision states $S_1$ and $S_2$, as
sketched in Figure \ref{chang-variable}. This curve is found by
introducing the linearization of the dynamical system around the
spontaneous equilibria $S_0$ as a new set of variables defined by
$X=P^{-1}(\nu-S_0)$, with $X=(x,y)^T$, or equivalently  $\nu = S_0
+ P X$ with $\nu=(\nu_1,\nu_2)^T$. Here, $P$ is the matrix
diagonalizing the jacobian of the dynamical system \eqref{ODE} at
the spontaneous state $S_0$, see \eqref{jac} in the Materials and
Methods section for more details, or \cite{CCM2}.

\begin{figure}[H]
\centering{\includegraphics[width=4in]{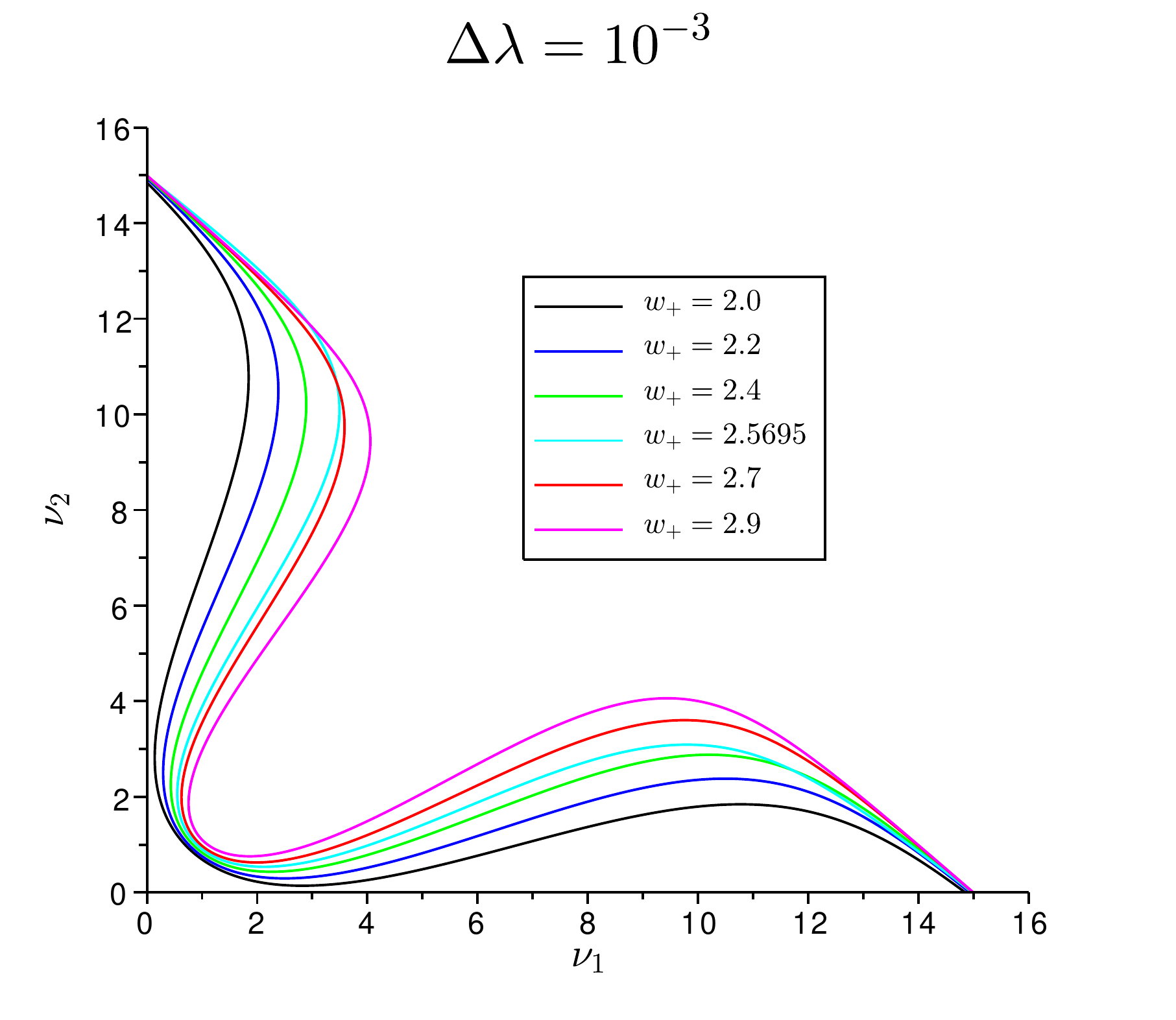}}
\caption{{\bf Slow manifold.} Variation of the approximated slow manifold represented on the $(\nu_1,\nu_2)$ plane for $w_+=2.0, \ldots , 2.9$ and for the biased case $\Delta \lambda=10^{-3}$. For $w_+ <1.9$ the slow manifold crosses the coordinate axes for small values of $\nu_1$ and $\nu_2$, becoming negative. The reduction theory is valid only for ranges of $w_+>1.9$.} \label{manifold}
\end{figure}

This change of variables, sketched in Figure \ref{chang-variable},
is a natural way to introduce slow $y$ and fast $x$ variables in
the system determined by the eigenvalues of the linearization at
the spontaneous state $S_0$. By rewriting the dynamical system
above \eqref{ODE} in these new variables together with the new
time variable $\tau=\varepsilon t$, we have an evolution governed by
\begin{equation*}%\label{ODExynew}
\begin{cases}
\varepsilon\displaystyle\frac{dx}{d\tau} = f(x,y) \\[4mm]
\displaystyle\frac{dy}{dt} = g(x,y)
\end{cases}\, ,
\end{equation*}
where $\epsilon$ is a small parameter. Assuming that the scaling
ratio is zero, $\epsilon=0$, we find the implicit relation
$f(x,y)=0$ that determines the approximated slow manifold
$x^*(y)$. In Figure \ref{manifold}, we plot the curves $x^*(y)$
for different values of $\omega_+$ with $\Delta \lambda=10^{-3}$.
We remark that this restricts the validity of this particular
example since these curves exit the set of positive values for the
rate variables. This comes from the fact that there are some
trajectories of the dynamical system without noise leading to
negative values for the rate variables. In this example the
approximated curves lie on the positive quadrant for $w_+\geq 1.9$
approximately.

\begin{figure}[H]
\centering{\includegraphics[width=3.1in]{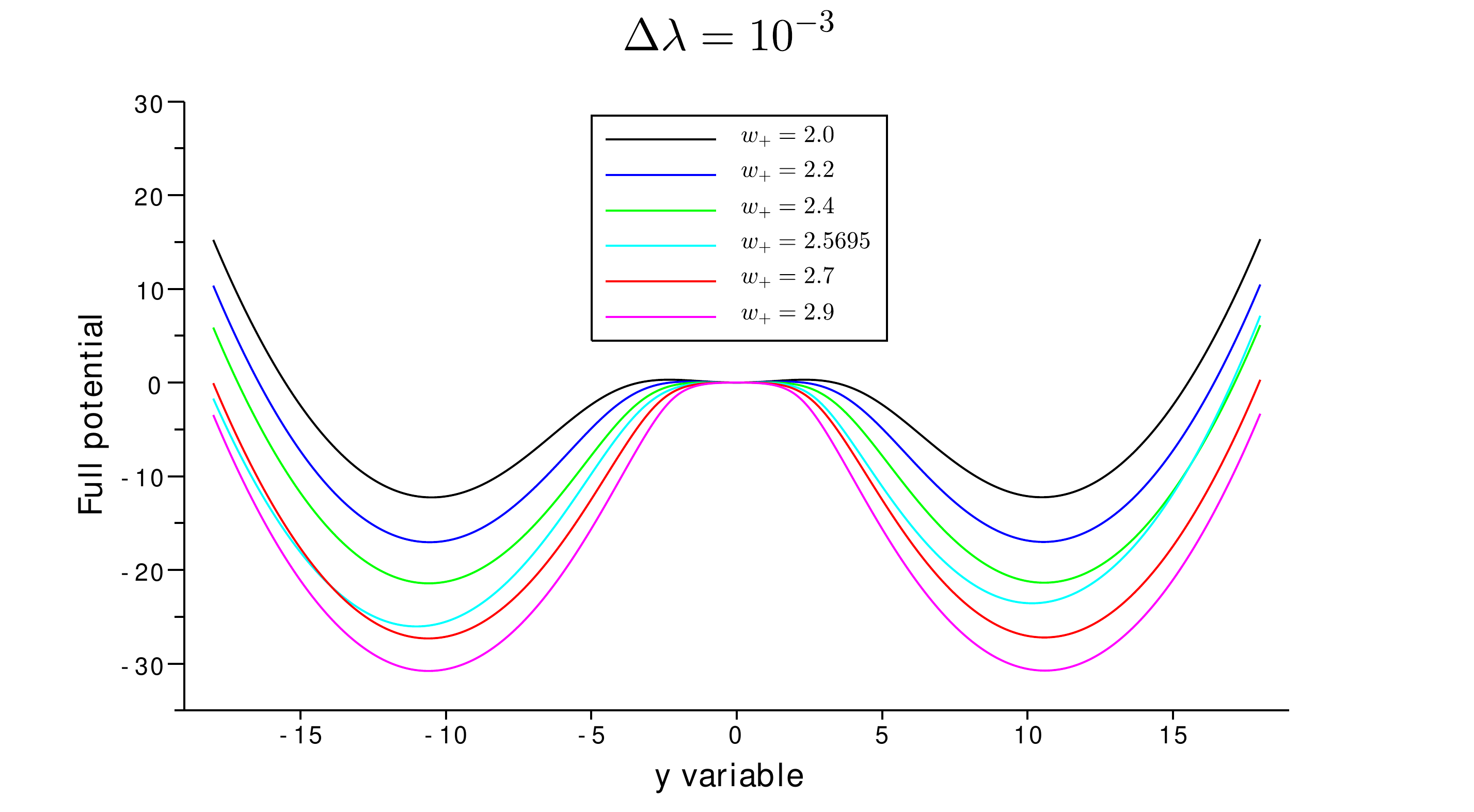}\;
\includegraphics[width=3.1in]{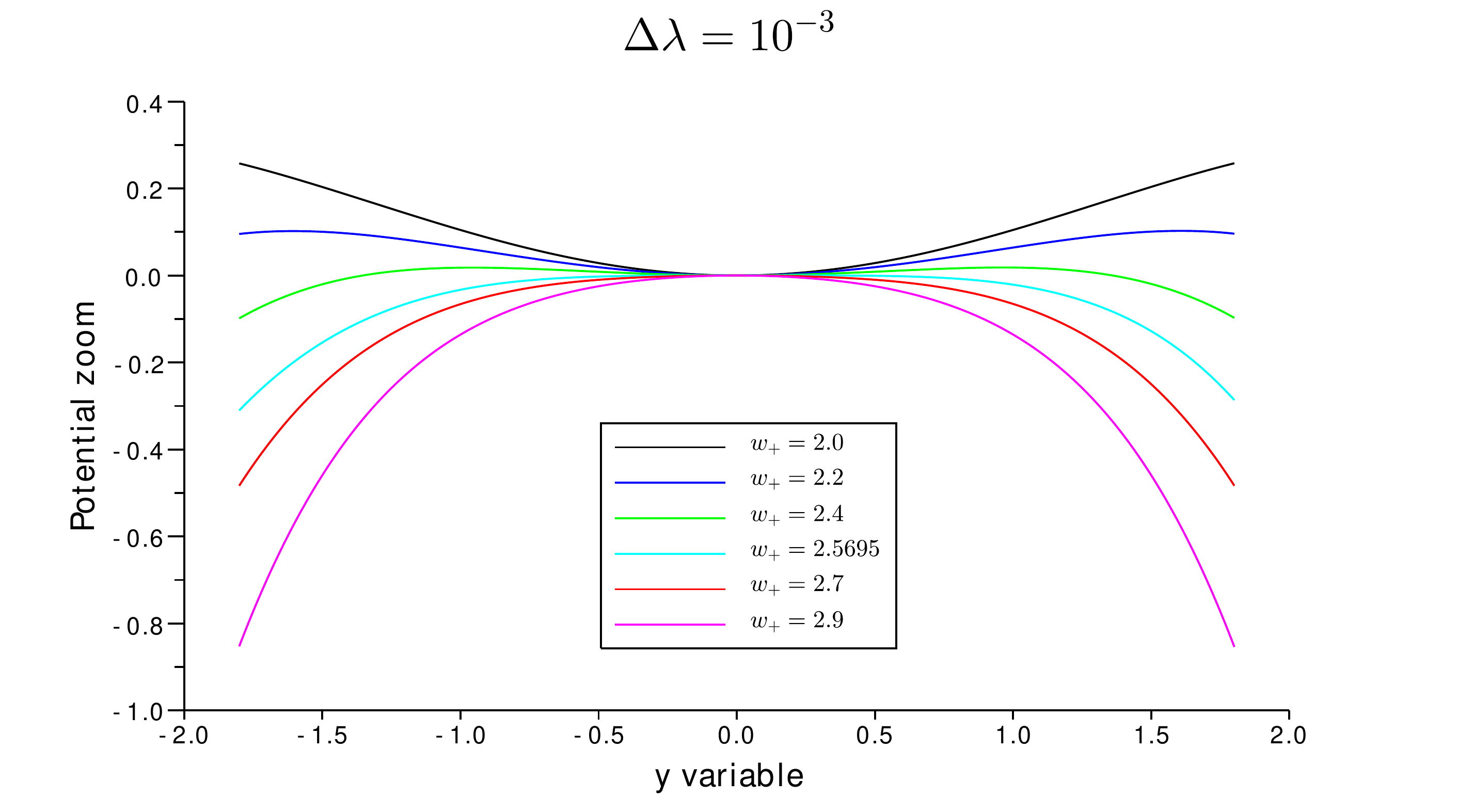}}
\caption{{\bf The 1D potential.} Variation of the approximated 1D potential traced along the $y$ variable, varying $w_+=2.0,
\ldots, 2.9$ and for the biased case $\Delta \lambda=10^{-3}$ (left). A zoom around the spontaneous state (right). Increasing $w_+$, the central equilibrium point passes from a minima to a maxima.}
\label{pot1D-all}
\end{figure}

Once we have found the slow manifold approximation $x^*(y)$, we can
restrict the dynamics in \eqref{ODE} to a single
effective Langevin equation. This equation is determined by a
potential $G(y)$ obtained from the evaluation of the dynamics over
the slow manifold approximation, leading to
$$
{dy} = g\left(x^*(y),y\right) {dt} + \beta_y \, dW_t^2\,,
$$
with $\beta_y$ properly obtained in terms of $\beta$ via the
change of variables. The effective potential is determined by the
relation $\partial_y G(y)=-g(x^*(y),y)$. The 1D effective computed
potential $G(y)$ for the biased case ($\Delta \lambda=10^{-3}$)
with respect to the slow variable $y$ is plotted in Figure
\ref{pot1D-all} for $w_+=2.0, \ldots, 2.9$ (left), together with a
zoom at the spontaneous state (right). We note that increasing
$w_+$ beyond the second bifurcation point, the spontaneous state
pass from local minimum to local maximum. We remind that the
effective potential $G(y)$ is given by formula \eqref{defG}
based on the approximated slow manifolds computed in Figure
\ref{manifold}.

\begin{figure}[!ht]
\centering{\includegraphics[width=3.1in]{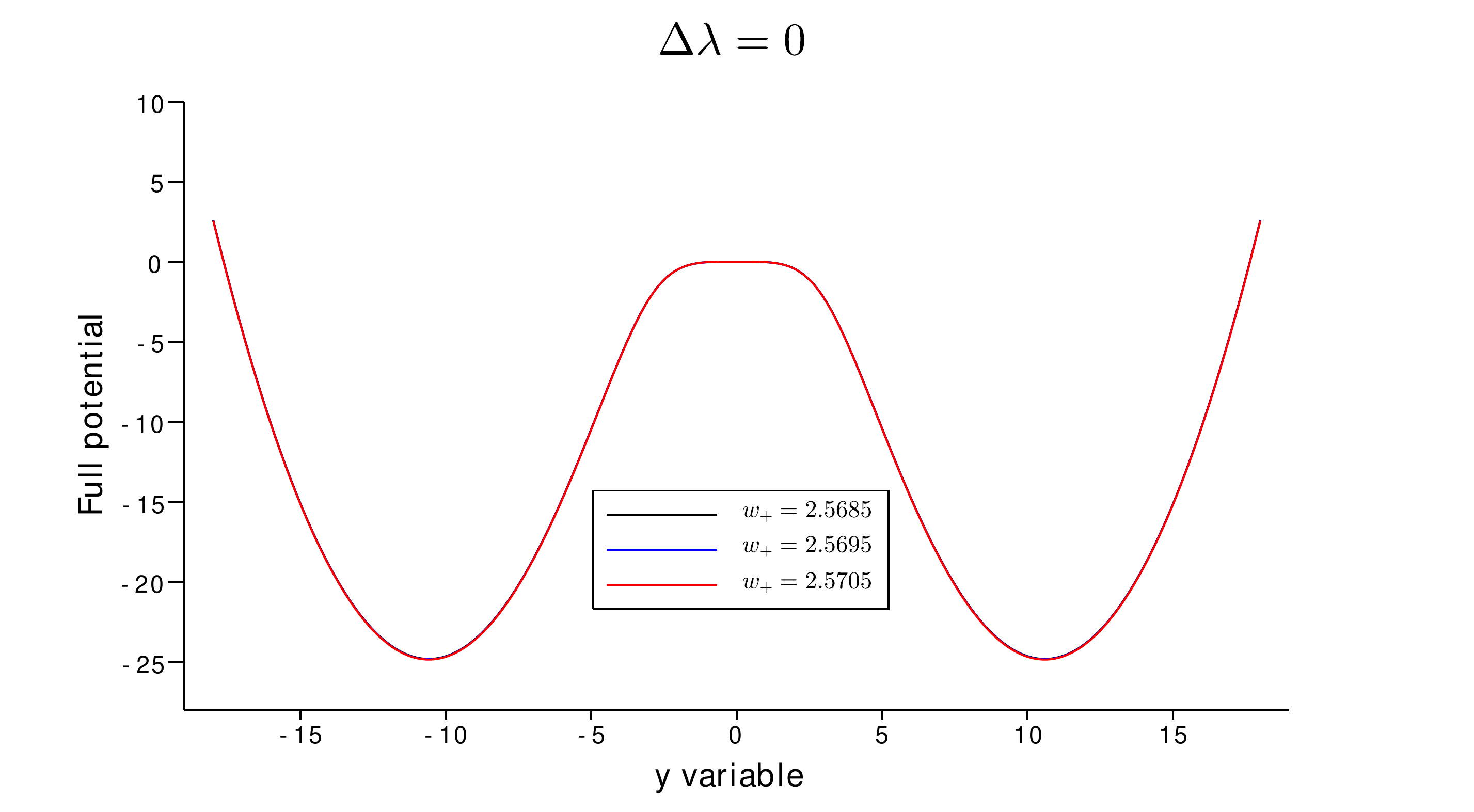}\;
\includegraphics[width=3.1in]{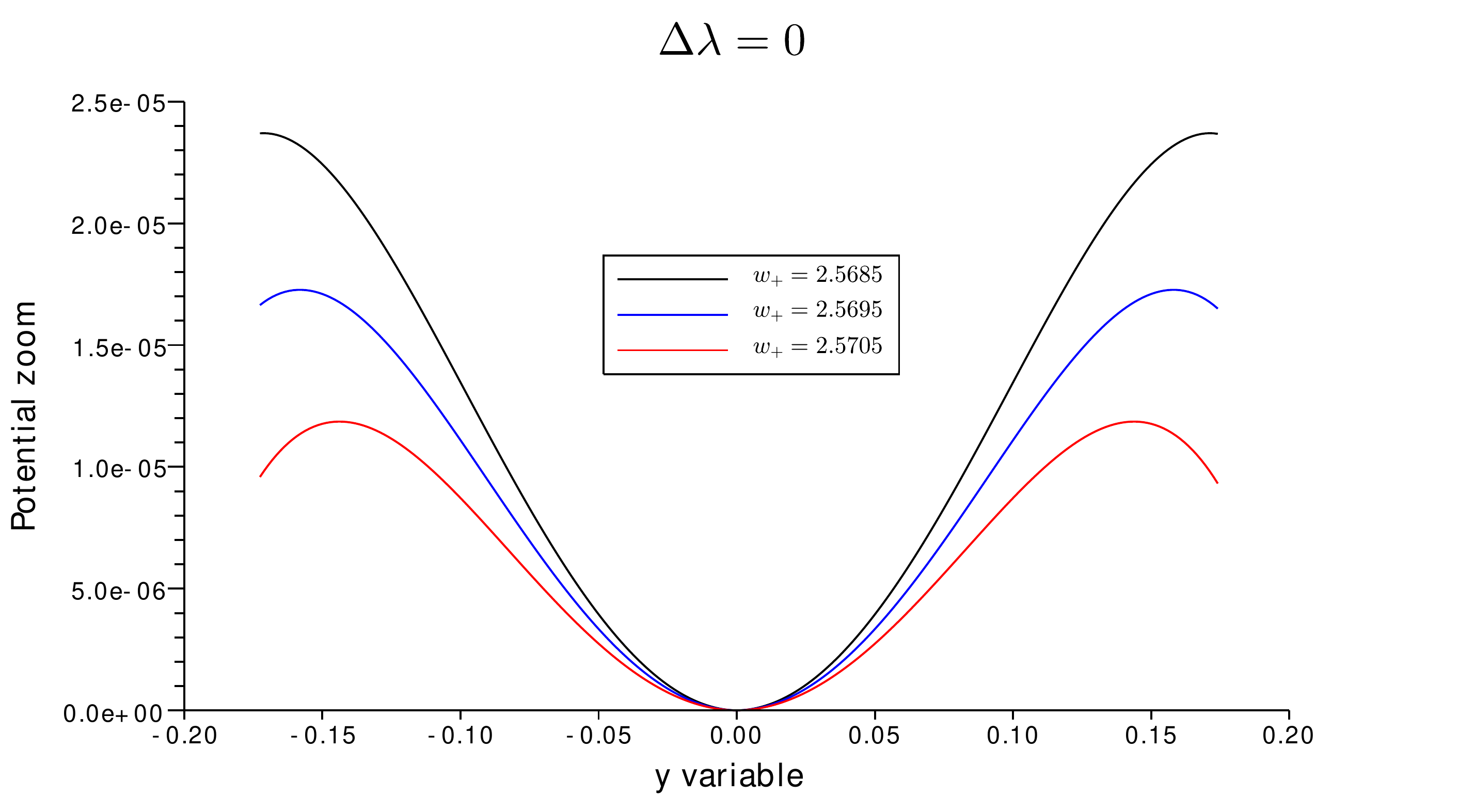}} \\
\centering{\includegraphics[width=3.1in]{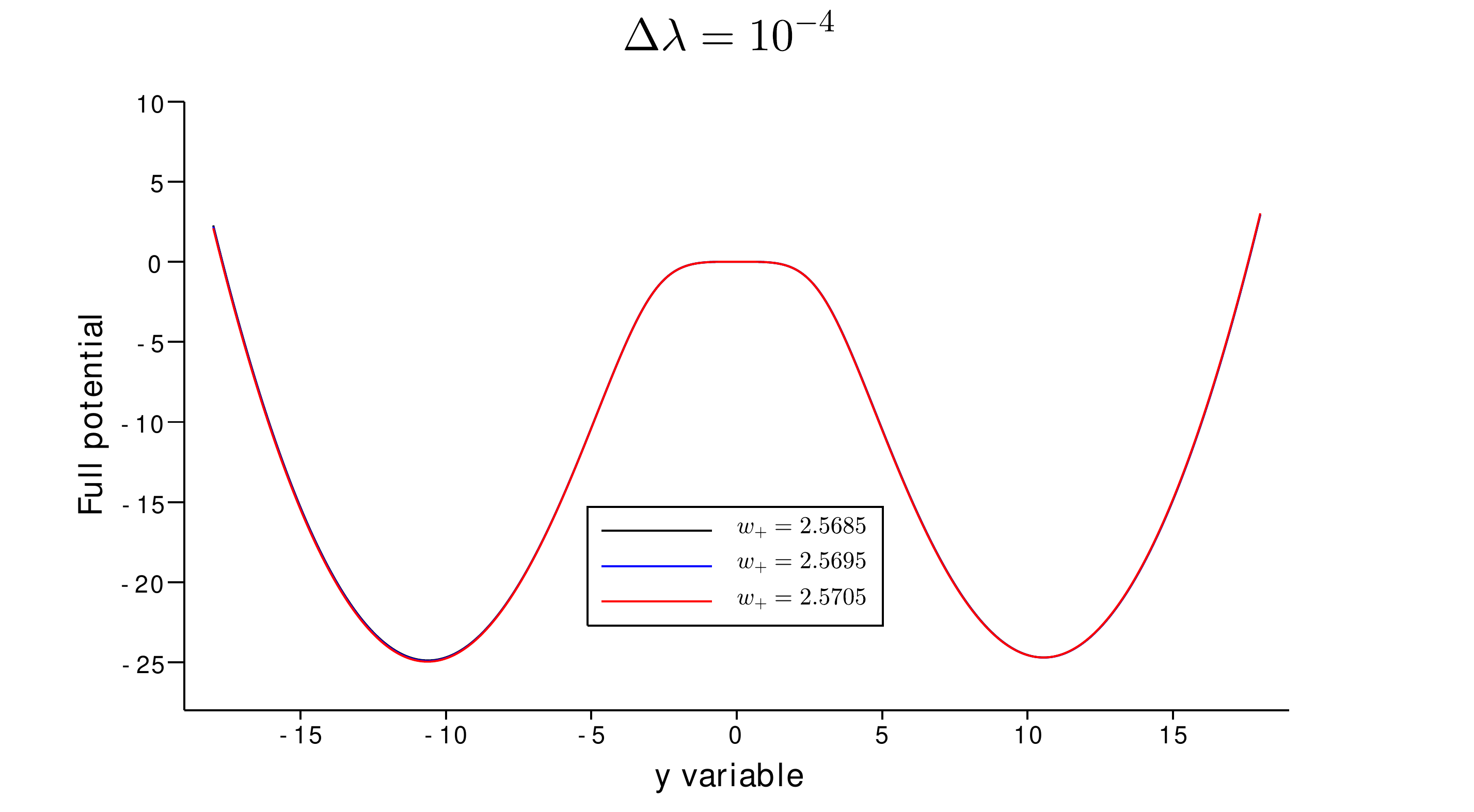}\;
\includegraphics[width=3.1in]{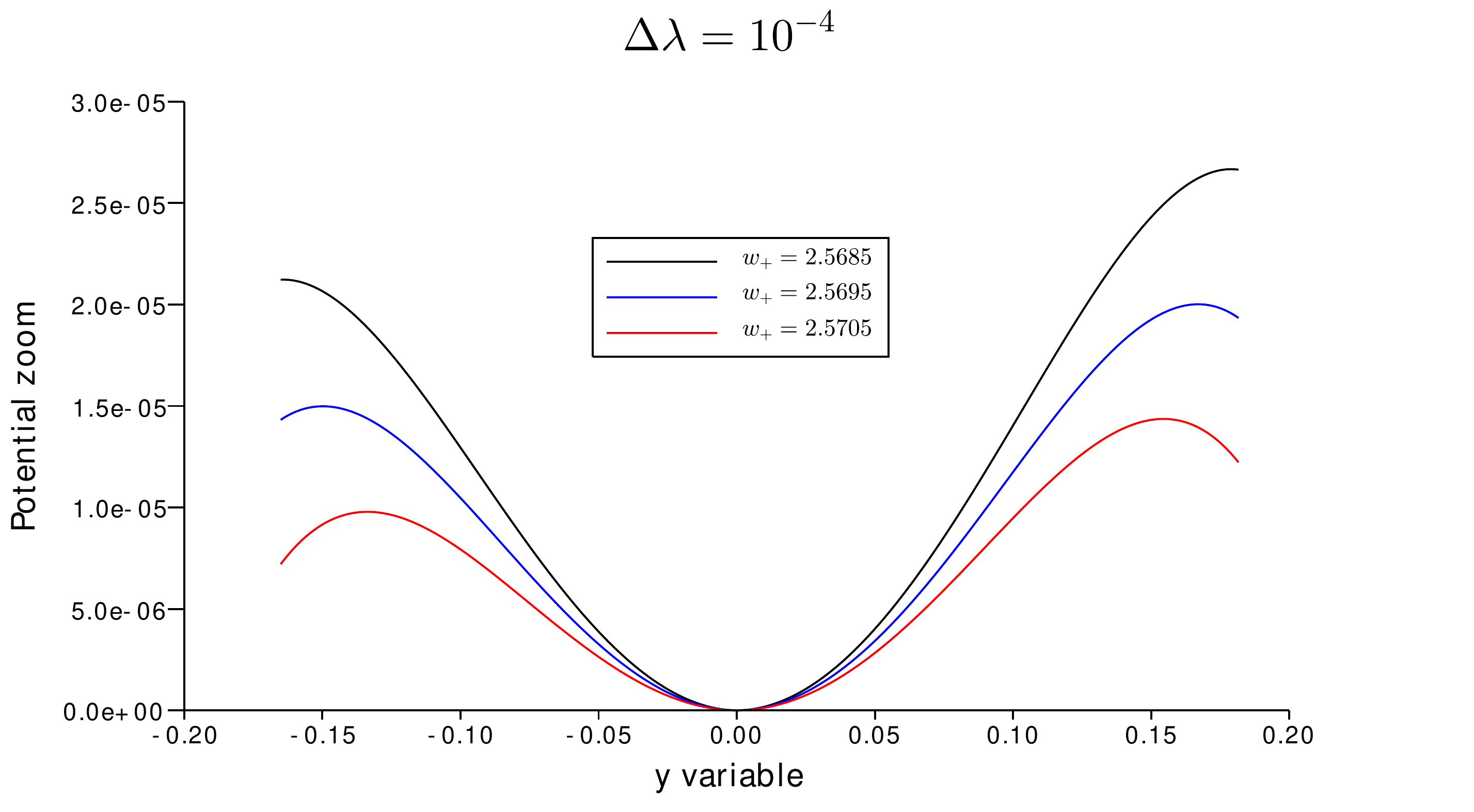}}\\
\centering{\includegraphics[width=3.1in]{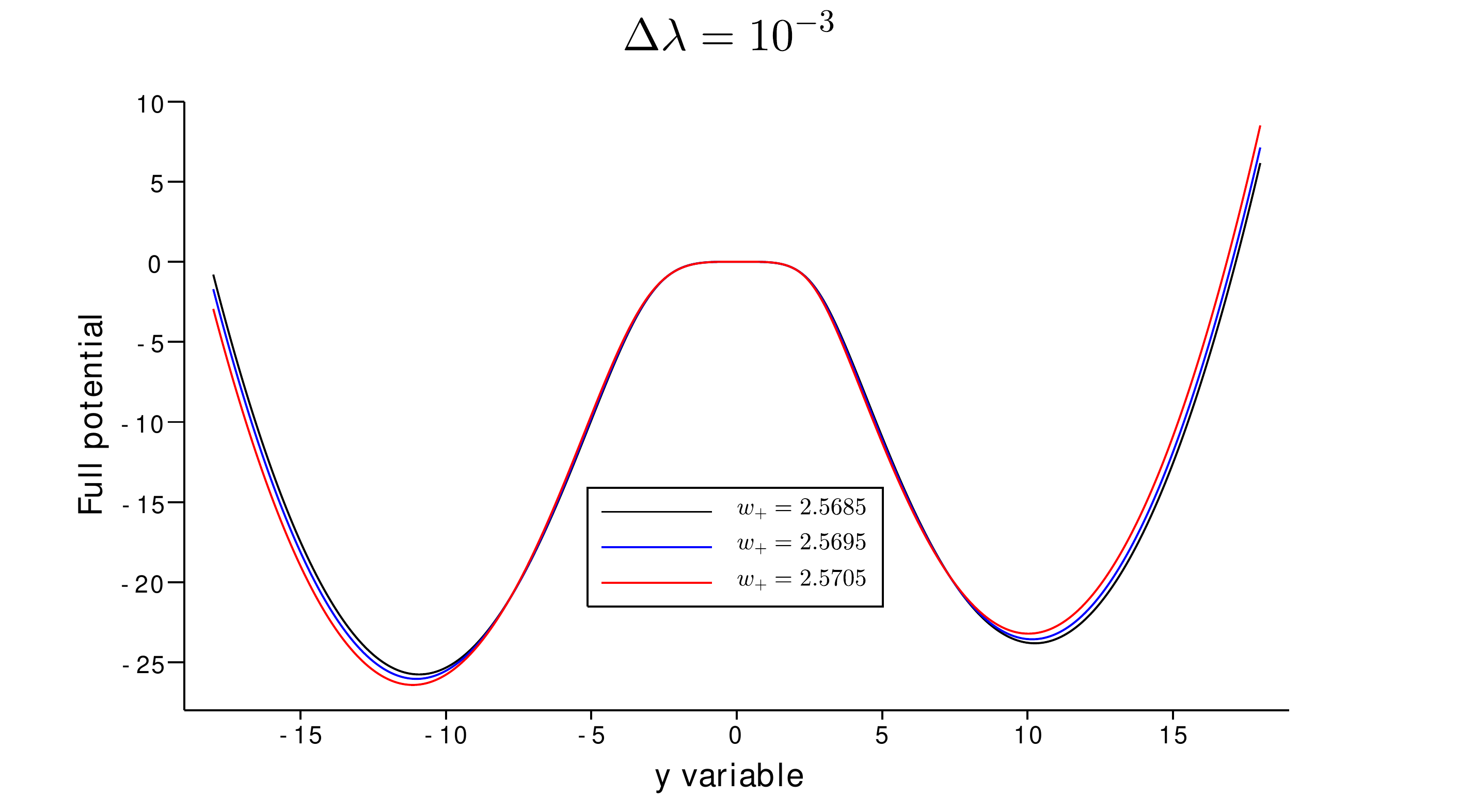}\;
\includegraphics[width=3.1in]{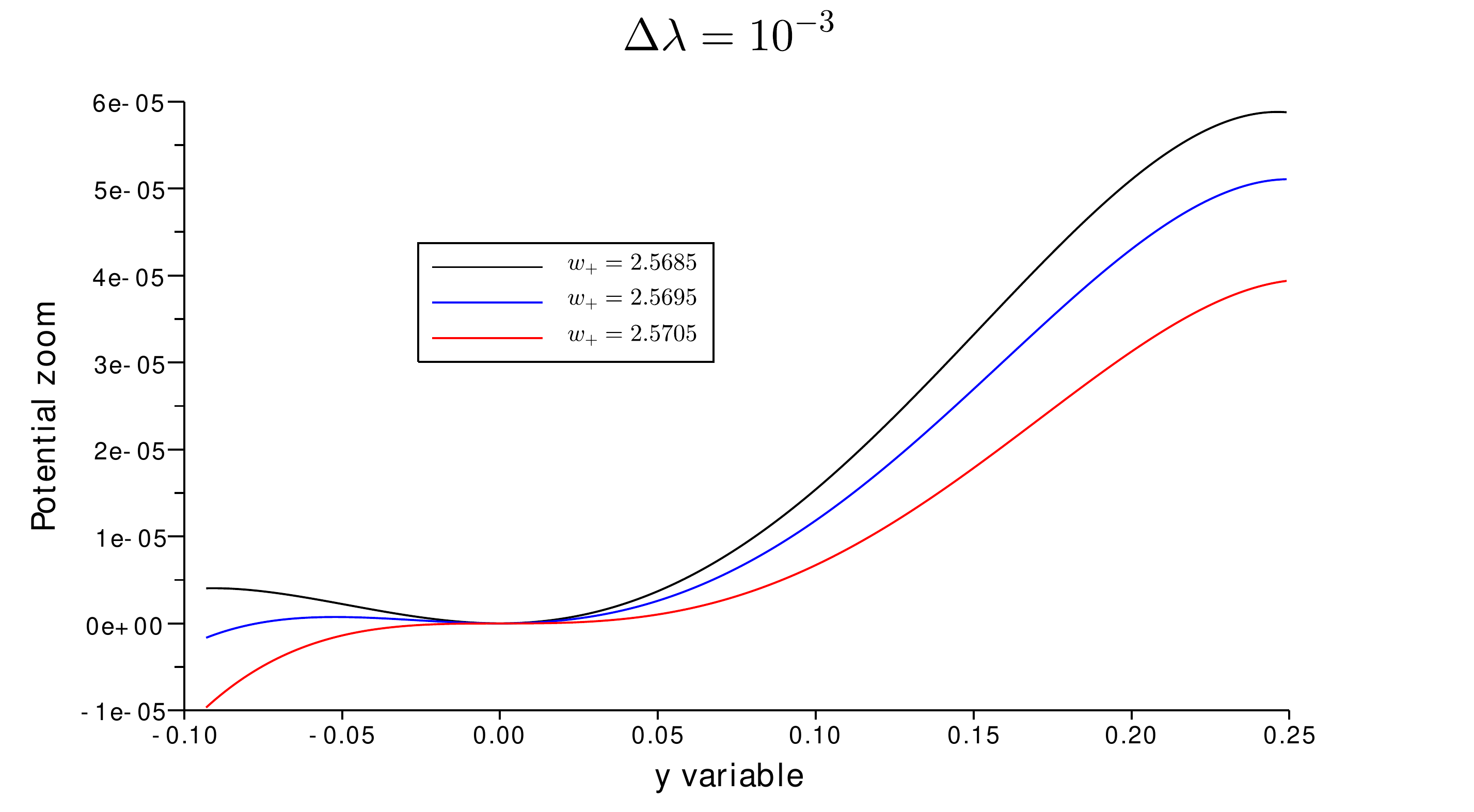}}
\caption{{\bf Potential at the bifurcation parameter $w_+=2.5696$}. Top left:
$\Delta \lambda=0$, top right: zoom for $\Delta \lambda=0$ around
the spontaneous state,  middle left: $\Delta \lambda=10^{-4}$, bottom
right: zoom for $\Delta \lambda=10^{-4}$ around the spontaneous
state, bottom left: $\Delta \lambda=10^{-3}$, bottom
right: zoom for $\Delta \lambda=10^{-3}$ around the spontaneous
state. For large bias $\Delta \lambda$ the potential becomes asymmetric, and for $\Delta \lambda =10^{-3}$ the spontaneous state changes its nature from a minima to a maxima.}
\label{pot1D}
\end{figure}

The 1D effective computed potentials $G(y)$ both for the unbiased
$\Delta \lambda=0$ and the biased case $\Delta \lambda = 10^{-4},
10^{-3}$ with respect to the slow variable $y$ are plotted in
Figure \ref{pot1D} for the three values of $w_+$ close to the
second bifurcation point ($w_+=2.5685$, $w_+=2.5695$,
$w_+=2.5705$). We note the asymmetry of the effective potential
$G(y)$ when increasing the bias $\Delta \lambda$ and the change in the nature of the spontaneous state when increasing the excitatory coefficient $w_+$.

This 1D effective potential can be considered as a good
approximation of the equilibrium state of the 2D Fokker-Planck
equation (associated to the system \eqref{ODE}) projected on the
1D-slow manifold which is parameterized by the slow variable $y$,
see Figure \ref{manifold}. We note that $y>0$ is pointing toward
$S_1$ while $y<0$ is pointing toward $S_2$, see Figure
\ref{chang-variable}. Therefore, the top left (respectively low
right) well whose minimum is attained at the stable decision state
$S_1$ (respectively $S_2$) in the 2D $\nu$-plane in Figure
\ref{chang-variable} corresponds to the right (respectively left)
well in Figure \ref{pot1D}.

\begin{figure}[!ht]
\centering{\includegraphics[width=4in]{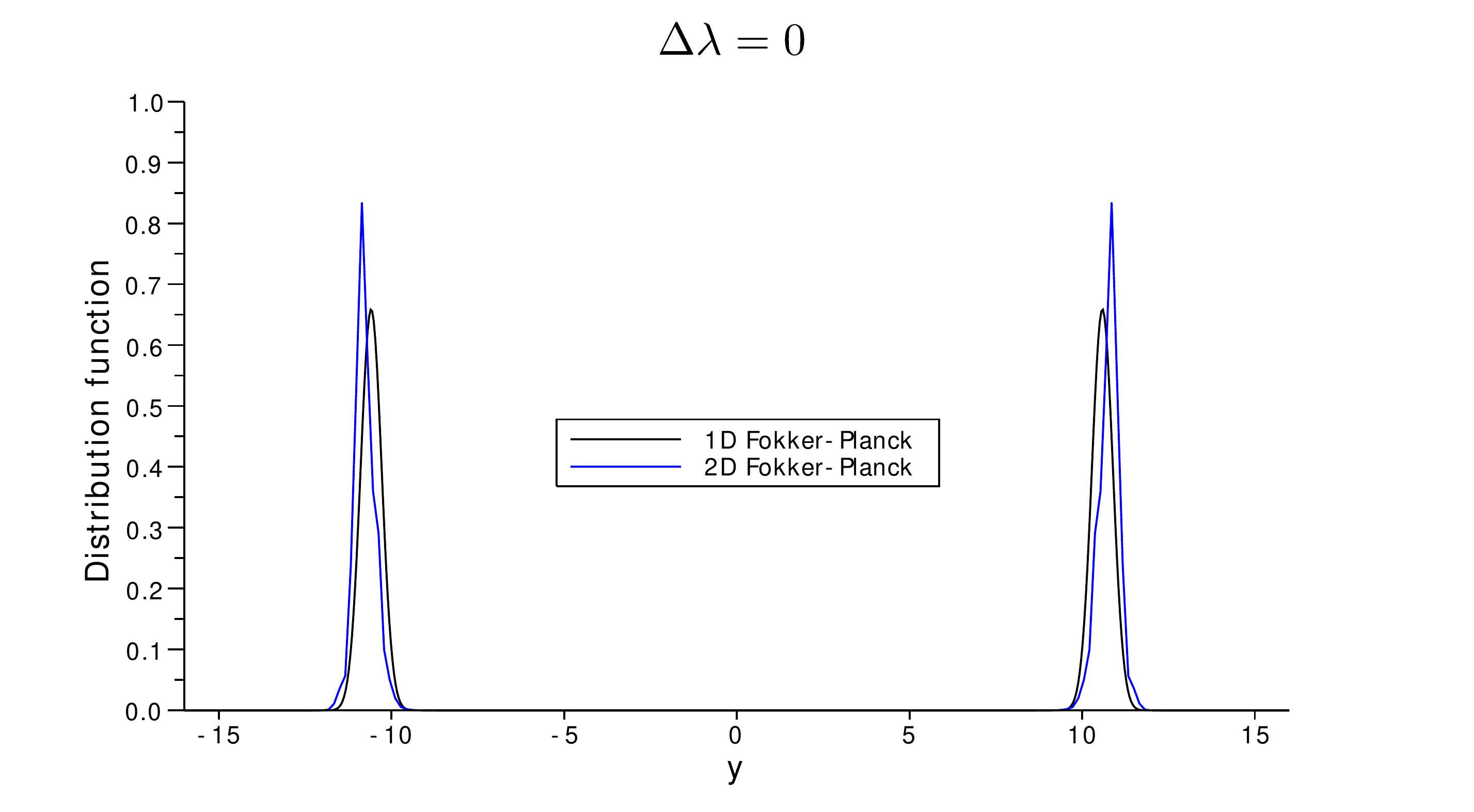}}\\
\centering{\includegraphics[width=4in]{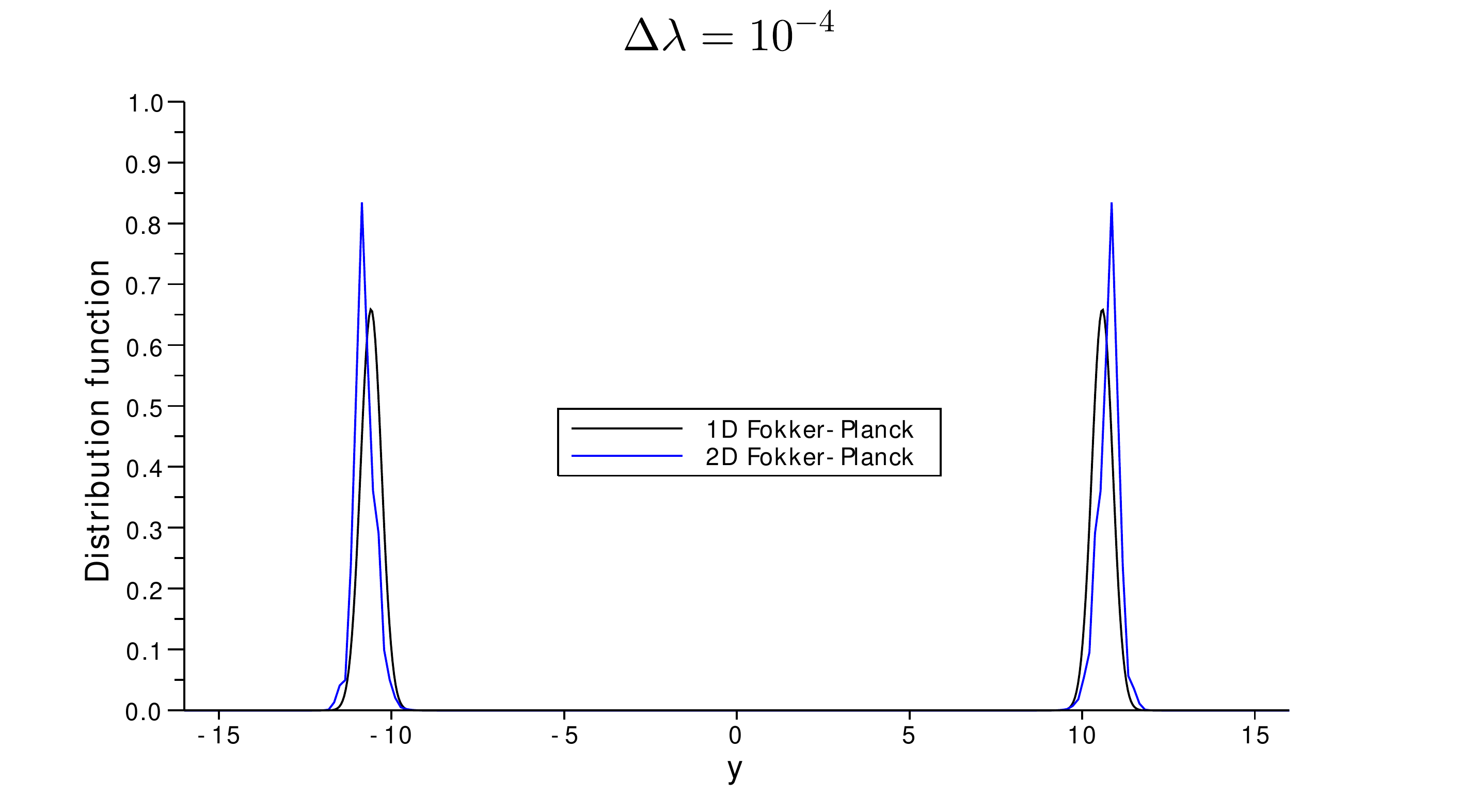}}\\
\centering{\includegraphics[width=4in]{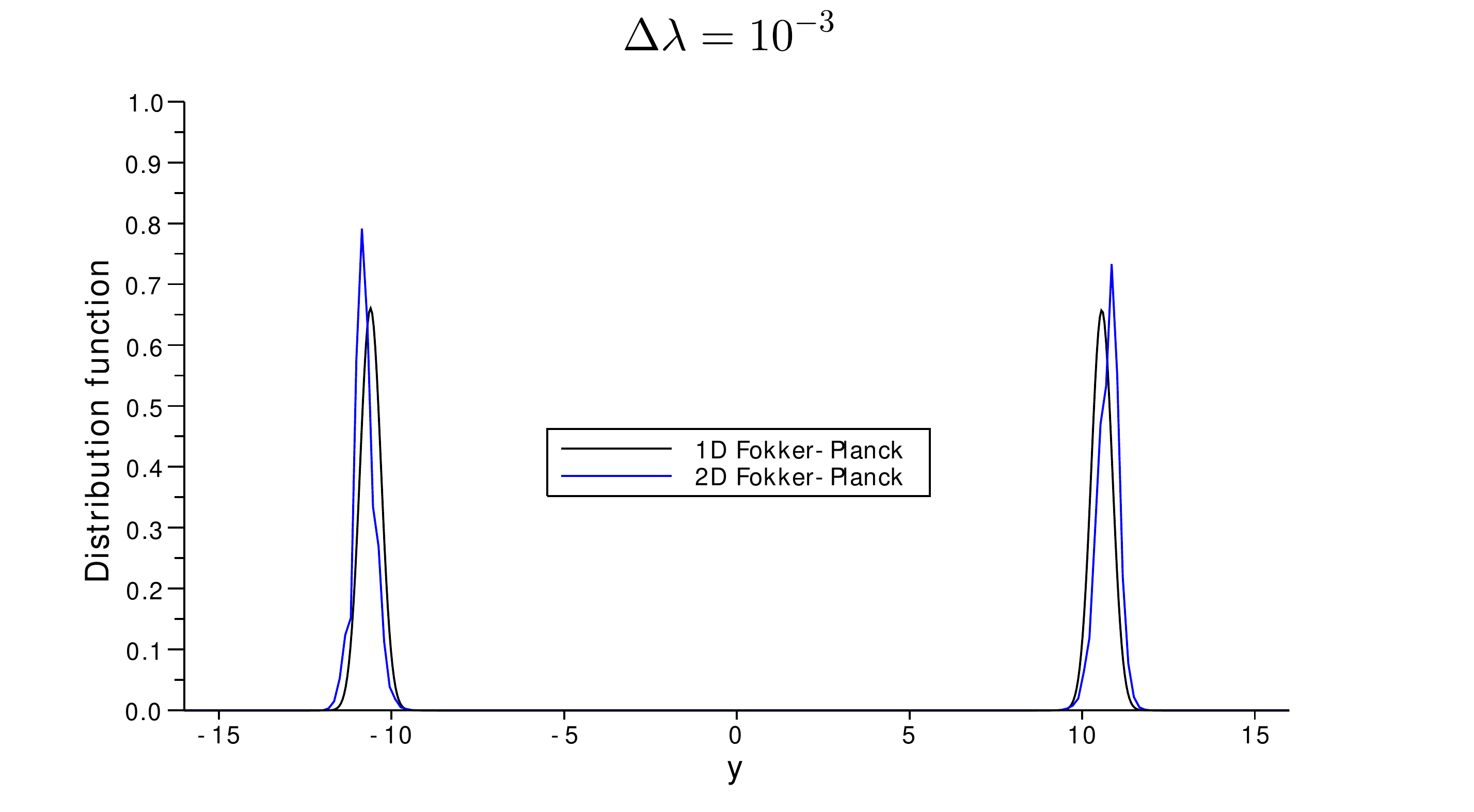}}
\caption{{\bf Distribution function comparison along $y$}. The computed marginals of the 2D problem (blue line) at time $t= 200 \ ms$ and the corresponding solution of the 1D reduced problem (black line) at the same time  show a good agreement. In this test we have chosen the bifurcation parameter $w_+=2.45$, the standard deviation $\beta=0.3$ and the bias $\Delta \lambda$ takes the values $0, 10^{-4}$ and  $10^{-3}$.}
\label{y-marginals2D}
\end{figure}

In Figure \ref{y-marginals2D}, we have computed the marginals
along $y$ of the 2D numerical computations of the Fokker-Planck
dynamics for $t= 200 \ ms$, compared with the solution of the 1D
reduction Fokker-Planck evolution equation \eqref{FPy} at $t= 200
\ ms$. We observe the errors committed by the 1D reductions and we
also observe that the slow-fast behavior is also present in the 1D
reduced Fokker-Planck equation. However, since the one dimensional
reduction is easily solved by numerical methods, then we can
achieve much larger computational times and compute explicitly the
approximated stationary state.

Concerning the comparison between the numerical solutions of the
2D complete Fokker-Planck and those of the 1D reduction, we refer
the reader to \cite{CCM2}. We just note here that since the
standard deviation is rather small ($\beta = 3 \cdot 10^{-3}$),
the stable states tends to be Dirac masses and the computational
costs are unbearable to compute 2D solutions approaching the
steady state. Note also that the 1D Fokker-Planck reduction can be
solved numerically applying an implicit in time scheme, drastically
reducing computational time costs.

\subsection*{Bifurcation: Reaction Time and Performance}

Finally, once the 1D reduction is validated by the above
arguments, the biologically interesting quantities to compute are
the performance $P_a$ and the reaction times $RT$. The first one
corresponds to the ratio of good answers at equilibrium, and the
second one to the smallest time needed to give an answer, no
matter wether it is good or wrong. The computations of the
reaction time $RT$ and performances $P_a$ use the formulas in the
supplementary material in \cite{RL}, which hold locally around the
spontaneous point and are based on the knowledge of the effective
potential $G(y)$.

\begin{figure}[!ht]
\centering{\includegraphics[width=4in]{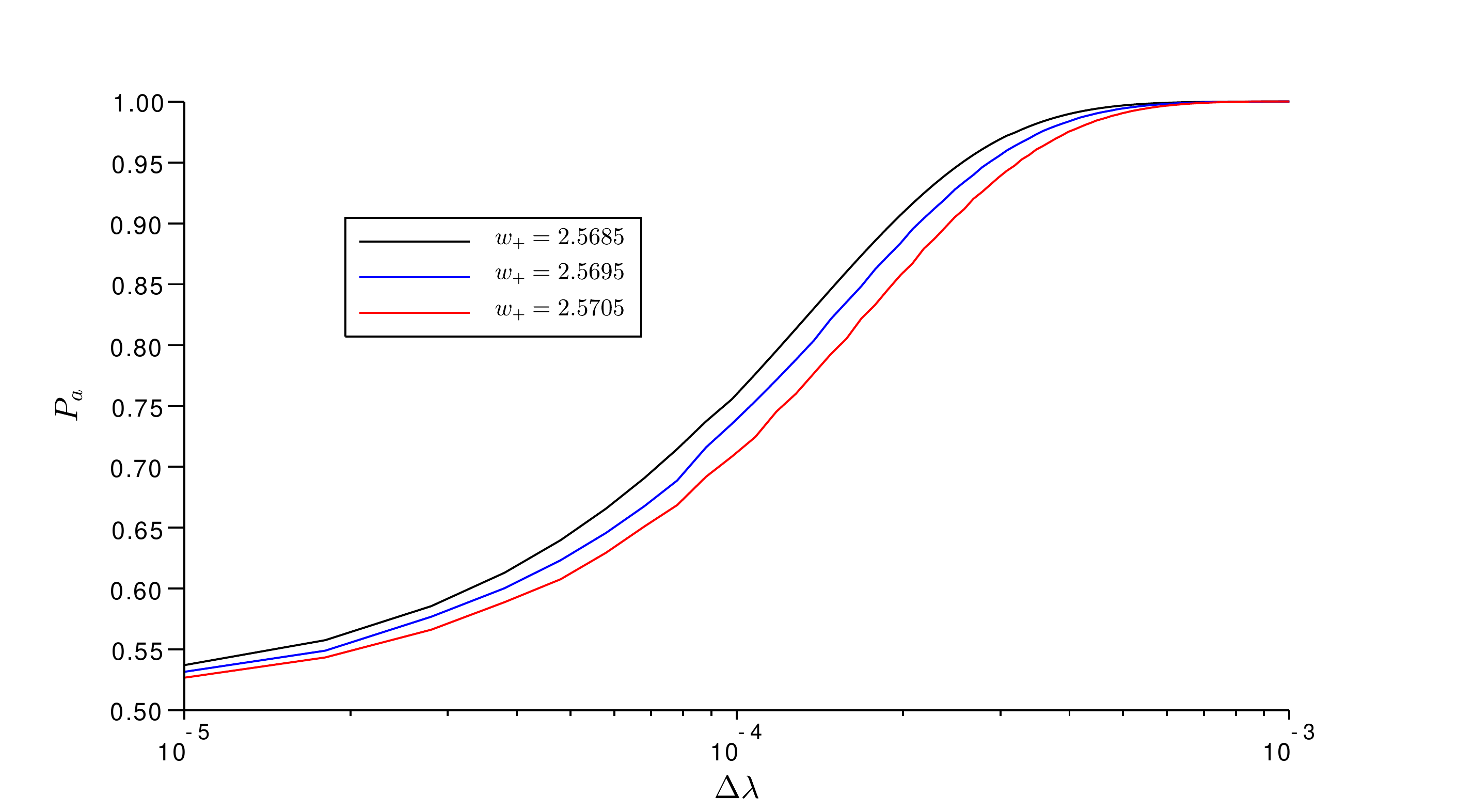}}\\
\centering{\includegraphics[width=4in]{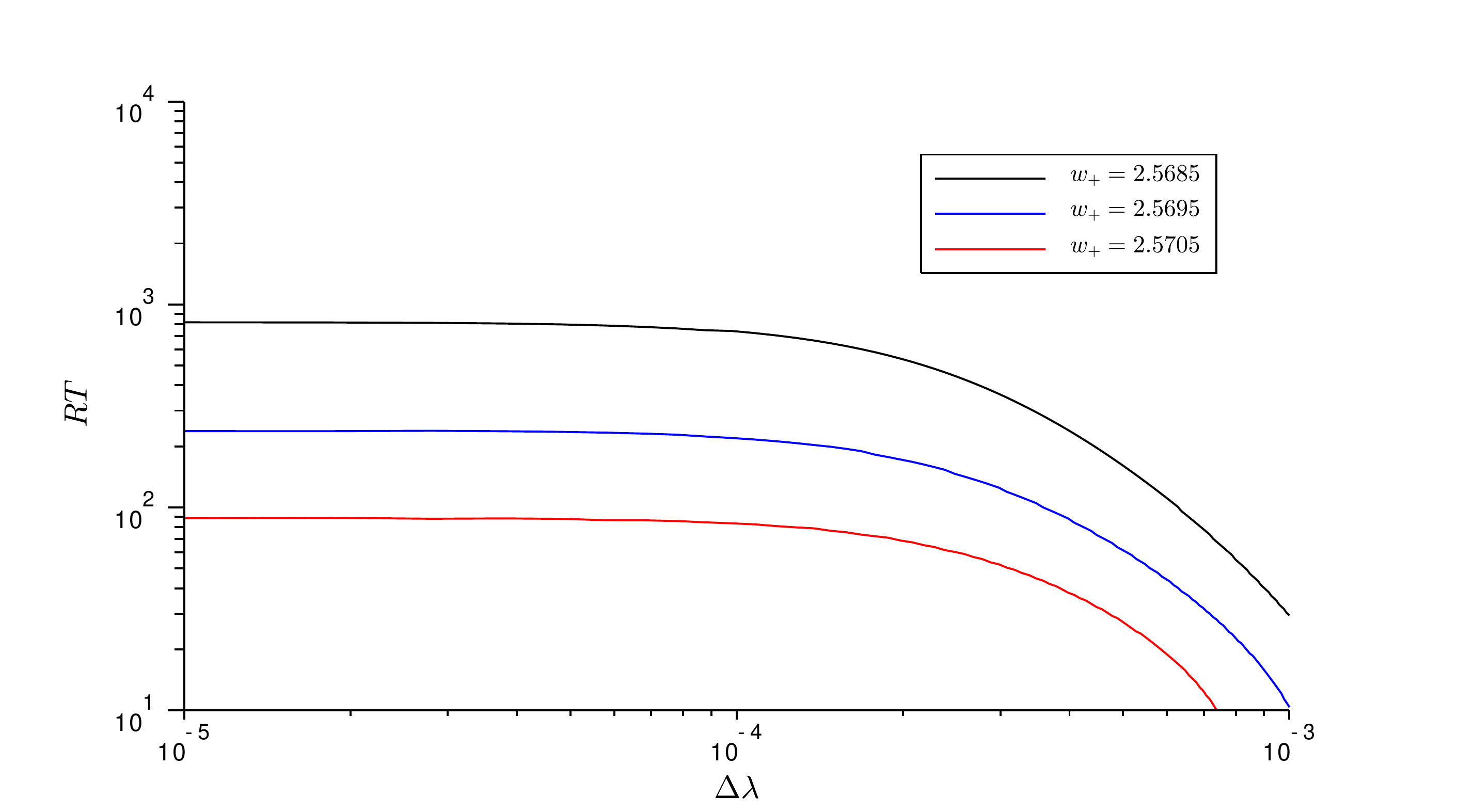}} 
\caption{{\bf Performance and Reaction Times.} Top: the performances ($P_a$) computed for three values of the bifurcation parameter $w_+$ close to the one of the spontaneous state and as a function of the bias $\Delta \lambda$. When $w_+$ increases, the performance decreases, and it converges to $1$ for large values of the bias $\Delta \lambda$. Bottom: the reaction times ($RT$) computed for three values of the bifurcation parameter $w_+$ close to the one of the spontaneous state and as a function of the bias $\Delta \lambda$. When $w_+$ increases, the reaction time decreases. In both pictures the standard deviation is $\beta=3\cdot 10^{-3}$. These behaviors are in agreement with previous results and experimental data.} \label{RT-PA}
\end{figure}

In Figure \ref{RT-PA} we plot the performance $P_a$ and reaction
time $RT$ for $\beta=3\cdot 10^{-3}$ as a functions of $\Delta
\lambda$ and $w_+$. They are computed for the decision state with
higher probability, the one for negative values of the
$y$-variable in Figure \ref{pot1D}, and over an interval
$[a^-,a^+]$ where $a^-$ and $a^+$ correspond to the
$y$-coordinates to the left and right of the spontaneous state for
which $G(y)$ has a maximum.  We observe the same qualitative
behavior as in the comparisons with experimental data performed in
\cite{RL}. As the bias increases, the potential gets tilted
towards the preferred decision state, and thus the performance of
the solution gets higher. Accordingly, the reaction time decreases
as the bias increases.

\section*{Discussion}

Here, we have reinforced the link between the underlying
physiology and the observed behavioral response in decision-making
tasks by developing a new strategy for the 1D reduction of the
dynamics of the neuronal circuit. In this way we have derived from
the underlying detailed neuronal dynamics a "nonlinear" diffusion
model valid for a wide range of the bifurcation parameter. We
observe how our method recovers the reduced one dimensional
dynamics proposed by Roxin and Ledberg in \cite{RL} close to the
bifurcation point since we keep the same qualitative behavior in
terms of performance and reaction times. We have also shown that
our reconstructed effective potential leads to a good
approximation of the stationary state of the two dimensional
dynamics projected on the slow manifold. This reduction allows for
efficient computations of the dynamics as soon as the approximated
slow manifold is well-defined.

% You may title this section "Methods" or "Models".
% "Models" is not a valid title for PLoS ONE authors. However, PLoS ONE
% authors may use "Analysis"
\section*{Materials and Methods}

In this section, we give more details of the 1D
reduction of system \eqref{ODE} presented above. Preliminary
results of this strategy were previously reported in \cite{CCM2}
in the easier case of the supercritical Hopf bifurcation. The
slow-fast behavior of the system (\ref{ODE}) with no noise, i.e.
$\beta=0$, can be characterized by the fact that the jacobian of
the linearized system at the unstable critical point
$S_0=(\nu_1^*,\nu_2^*)^T$ has a large condition number. More
precisely, we write the deterministic part of the dynamical system
\eqref{ODE} as $\dot \nu = F(\nu)$, where $\nu=(\nu_1,\nu_2)^T$ is
the rate vector and $F(\nu)=-\nu+\Phi(\Lambda + W\nu)$ is the
drift. The linearized jacobian matrix at any equilibrium point is
given by:
\begin{equation}\label{jac}
J_F(z_1,z_2)=\left(
\begin{array}{cc}
-1 + w_{11} \varphi'(z_1)  &  w_{12} \varphi'(z_1) \\
w_{21}  \varphi'(z_2) & -1 + w_{22} \varphi'(z_2)
\end{array}
\right) ,
\end{equation}
where we have denoted by $z_i$ the values $z_i := \lambda_i +
w_{i1} \nu_1 + w_{i2} \nu_2$. Since $S_0$ is a hyperbolic fixed
point (saddle point), the jacobian $J_F(S_0)$ has two real
eigenvalues, $\mu_1$ and $\mu_2$ of opposite signs. Let us denote
by $\mu_1$ the (large in magnitude) negative eigenvalue and by
$\mu_2$ the (small) positive eigenvalue of $J_F(S_0)$. The small
parameter $\varepsilon <<1$ responsible for the slow-fast behavior
is determined by the ratio
\begin{equation*}%\label{e}
\varepsilon= -\mu_2/\mu_1.
\end{equation*}
In order to reduce the system we need to introduce new variables
based on the linearization of the problem. We will denote by $P$
the matrix containing the normalized eigenvectors of $J_F(S_0)$
and by $\PP$ its inverse matrix. Furthermore, using the
Hartman-Gro$\beta$man theorem \cite{Ha}, we know that the
solutions of the dynamical system are topologically conjugate with
its linearization in the neighborhood of a hyperbolic fixed
point. Let us write it as follows
\begin{equation}\label{JF}
 J_F(S_0) = P D P^{-1},
\end{equation}
where  $D=diag(\mu_1,\mu_2)$
is the associated diagonal matrix. We can describe the coordinates
$\nu$ in the eigenvector basis and centered on the saddle point
$S_0$ as follows:
\begin{equation*}%\label{chvar}
\nu = S_0 + P X ,
\end{equation*}
which gives the definition for the new variable $X=(x,y)^T$, see
also figure \ref{chang-variable}, $X=P^{-1}(\nu-S_0)$.

In these new coordinates $x$ corresponds to the fast scale while
$y$ is the slow varying variable. We can conclude that system
\eqref{ODE} reads in the $X$ phase space as $\dot X = H(X)$ where
$H(X)=(f(x,y),g(x,y))$ is the two dimensional vector defined by
$H(X)=P^{-1} F(S_0 +P X)$, and thus,
\begin{equation}\label{ODExy}
\begin{cases}
\displaystyle\frac{dx}{dt} = f(x,y) := -x-(P^{-1}S_0)_1 + \left(P^{-1}
\Phi \left( \Lambda + W (S_0 + P X) \right)\right)_1\\[4mm]
\displaystyle\frac{dy}{dt} = g(x,y) = -y-(P^{-1}S_0)_2 + \left(P^{-1}
\Phi \left( \Lambda + W (S_0 + P X ) \right)\right)_2
\end{cases}\, .
\end{equation}
Direct computations gives that $ J_H(X) = P^{-1} J_F(S_0 +PX) P$,
and using \eqref{JF} and that $X(S_0)=0$, we obtain that $J_H(0) =
D$. We can choose a new time scale for the fast variable
$\tau=\varepsilon t$ in such a way that for $t\simeq 1$, then
$\tau\simeq \varepsilon^{-1}$. Then, the fast character of the
variable $x$ is clarified and the system reads as
\begin{equation*}%\label{tildeODExy}
\begin{cases}
\varepsilon\displaystyle\frac{dx}{d\tau} = f(x,y) \\[4mm]
\displaystyle\frac{dy}{dt} = g(x,y)
\end{cases} .
\end{equation*}
Our claim is that the curve of equation $f(x,y) =0$ approximates
well when $\epsilon \ll 1$ the slow manifold, that is the unstable
manifold that joins the spontaneous point $S_0$ to the other
equilibrium points. We refer the reader to
the discussion in \cite{BG} on stochastic slow-fast dynamics. Due
to the non-linearity of the function $f$ in \eqref{ODExy}, we
cannot expect an explicit formula for the solution of this
equation. Nevertheless, since $\partial_x f(0,0) \neq 0$, the
resolution in the neighborhood of the unstable equilibrium point
$(0,0)$ is insured by the implicit function theorem. Hence we can
define a curve:
\begin{equation}\label{def_curve}
x=x^*(y),
\end{equation}
such that $f\left(x^*(y),y\right)=0$ in a neighborhood of $(0,0)$. We note
also that, by construction the approximated slow manifold $S_0 + P
\left(x^*(y),y\right)^T$, implicitly defined by (\ref{def_curve}), intersects
the exact slow manifold at all equilibrium points, i.e. where both
$f$ and $g$ vanishes. Finally, we can conclude the slow-fast
ansatz, replacing the complete dynamics by the reduced dynamics on
the approximated slow manifold, and obtain the reduced one
dimensional differential equation:
\begin{equation*}%\label{slowODE}
\dot y = g\left(x^*(y),y\right) .
\end{equation*}

A similar treatment can be done in the presence of the noise
terms. By changing variables from $\nu$ to $X$ and since the new
variables $x$ and $y$ are linear combination of $\nu_1$ and
$\nu_2$, then it is standard to check
\begin{equation*}%\label{ODExynoise}
\begin{cases}
{dx} = f(x,y) {dt} + \beta \sqrt{ a_{11}^2 + a_{12}^2} \,dW_t^1\\[3mm]
{dy} = g(x,y) {dt} + \beta \sqrt{ a_{21}^2 + a_{22}^2} \, dW_t^2
\end{cases}\, ,
\end{equation*}
where $(a_{ij})$ are the elements of the matrix $\PP$ and $W_t^i$ are two independent normalized white noise. Arguing
as in the deterministic case, we can choose a fast time scale for
the variable $x$ and we note that the noise term also changes since
$\sqrt{\varepsilon} dW_t^1 = dW_\tau^1$ and $\varepsilon
dt=d\tau$. Therefore, we can deduce again that the reduced one
dimensional model must read:
\begin{equation}\label{ODEy}
{dy} = g\left(x^*(y),y\right) {dt} + \beta_y \, dW_t^2 ,
\end{equation}
with $\beta_y=\beta \sqrt{ a_{21}^2 + a_{22}^2}$.

Finally, we can consider the Fokker-Planck or forward Kolmogorov
equation associated to the 1D stochastic differential
equation \eqref{ODEy}. This gives the reduced dynamics for the
probability density $q(t,y)$ of finding neurons with rate
determined by the rate $y\in [-y_m,+y_m]$ over the approximated
slow manifold $S_0+P \left( x^*(y),y\right)^T$, for $t\geq 0$. Actually, it
must obey to the following 1D Fokker-Planck equation:
\begin{equation}\label{FPy}
\partial_t q + \partial_y \left(  g\left( x^*(y),y\right) q - \frac12 \beta_y^2 \partial_y q \right)  =0 ,
\end{equation}
with no-flux boundary conditions on $y=\pm y_m$: $g\left( x^*(y),y\right) q -
\frac12 \beta_y^2 \partial_y q =0$, see \cite{Gar} and \cite{CCM2}. The 1D
Fokker-Planck dynamics \eqref{FPy} are given for large times by
the stationary solution given by the invariant measure for the
stochastic differential equation \eqref{ODEy}. We can easily find
it by defining the associated potential $G(y)$ being the
antiderivative of the flux term $g\left( x^*(y),y\right)$. In other words, we
can always define the potential function:
\begin{equation}\label{defG}
G(y) = -\int_0^y g\left( x^*(z),z\right) dz .
\end{equation}
Then, the stationary probability density of the 1D Fokker-Planck
dynamics \eqref{FPy} is given by
\begin{equation*}%\label{qstat}
q_s(y)=\frac1{Z}\exp\left(-2G(y)/\beta_y^2\right)\, ,
\end{equation*}
where $Z$ is the normalization constant. As explained also in
\cite{CCM}, this stationary solutions are the asymptotic
equilibrium states for the solution of the Fokker-Planck equation.
In other words, letting time going to infinity, the solution
$q(t,y)$ to \eqref{FPy} must converge to $q_s(y)$. We have shown
in \cite{CCM} that the decay to equilibrium for the two
dimensional problem was exponential. This rate of convergence is
also true for the 1D reduction since the potential is a small
perturbation of convex potentials, see the well-know results in
\cite{HS,AMTU}. We also recall that we are interested in the long
time behavior of the solutions and that the convergence in the
slow manifold given by the variable $y$ is slow. Hence it is
relevant to have a direct computation of their asymptotic behavior
without need to solve the whole 1D or 2D Fokker-Planck equation.

We also remind that the solution of the system \eqref{ODE} is
related to the 2D Fokker-Planck or forward Kolmogorov equation:
\begin{equation}\label{FP}
\partial_t p + \nabla \cdot \left( F\, p  -  \bb \nabla p \right) =0
\end{equation}
with boundary conditions $\left( F\, p - \bb \nabla p \right)
\cdot n =0$ on the boundary of the domain $[0,\nu_m]^2$, see
\cite{Gar}. Let us remark that, once we know the stationary state
of the 1D reduced Fokker-Planck equation \eqref{FPy}, we can
approximate the long-time dynamics of all quantities of interest
related to the system \eqref{ODE} and the 2D Fokker-Planck
\eqref{FP}. More precisely, as $\varepsilon\to 0$, $p$ approaches
a concentrated density along the the curve $\nu(y)=S_0+P(x^*(y),
y)^T$. Then, for any test function $\Psi=\Psi((\nu_1,\nu_2)^T)$,
the moment $M_\Psi$ of the stationary probability distribution
function $p_s(\nu_1,\nu_2)$ as $\varepsilon\to 0$ can be
approximated by
$$
M_\Psi = \int_\Omega \Psi(\nu) \, p_s \, d\nu_1 d\nu_2 \simeq \int
\Psi\left( S_0+P\left( x^*(y),y\right)^T\right) \, q_s(y) \,dy \, .
$$
This formulae can be used to compute either classical moments of
$p_s$ or marginals.

% Do NOT remove this, even if you are not including acknowledgments
\section*{Acknowledgments}
JAC acknowledges support from the Royal Society by a Wolfson
Research Merit Award and by the Engineering and Physical Sciences
Research Council grant with references EP/K008404/1. JAC was
partially supported by the project MTM2011-27739-C04-02 DGI
(Spain) and 2009-SGR-345 from AGAUR-Generalitat de
Catalunya. SC and SM acknowledge support by the ANR
project MANDy, Mathematical Analysis of Neuronal Dynamics,
ANR-09-BLAN-0008-01.

%\section*{References}
% The bibtex filename
\bibliography{biblio}

\end{document}